\documentclass[a4paper, 12pt]{article}

\usepackage{amsfonts, amsmath, amsthm}
\usepackage{mathrsfs, latexsym}

\allowdisplaybreaks

\theoremstyle{plain}
\newtheorem{thm}{Theorem}[section]

\newtheorem{lem}{Lemma}[section]

\theoremstyle{definition}

\theoremstyle{remark}

\theoremstyle{remark}

\numberwithin{equation}{section}

\newcommand{\R}{\mathbb{R}}
\newcommand{\C}{\mathbb{C}}

\newcommand{\F}{\mathcal F}
\newcommand{\IF}{\mathcal F^{-1}}

\newcommand{\M}{\mathcal M}

\newcommand{\D}{\mathcal D}

\newcommand{\B}{\mathcal B}

\newcommand{\W}{\mathcal W}

\newcommand{\G}{\mathcal G}
\newcommand{\IG}{{\mathcal G}^{-1}}

\newcommand{\U}{\mathcal U}

\newcommand{\pa}{\partial}
\newcommand{\eps}{\varepsilon}
\newcommand{\jb}[1]{\left\langle #1 \right\rangle}
\newcommand{\conj}[1]{\overline{#1}}

\DeclareMathOperator{\im}{\rm Im}

\begin{document}
\title{
A note on decay rates of solutions to a system of cubic nonlinear 
Schr\"odinger equations in one space dimension
}

\author{
          Donghyun Kim\thanks{
              Seoul National Science Museum. 
              215, Changgyeonggung-ro, Jongno-gu, Seoul 110-360, Korea
              (E-mail: {\tt u553252d@alumni.osaka-u.ac.jp})
            }
}

\date{\today}   
\maketitle

{\noindent{\bf Abstract:}\ 
We consider the initial value problem for a system of cubic nonlinear 
Schr\"odinger equations with different masses in one space dimension.
Under a suitable structural condition on the nonlinearity, 
we will show that the small amplitude solution exists globally 
and decays of the rate $O(t^{-1/2}(\log t)^{-1/2})$ in $L^\infty$ as
$t$ tends to infinity, if the system satisfies certain mass relations.\\}




\section{Introduction} \label{sec_intro}
We consider the initial value problem for a system of cubic nonlinear
Schr\"odinger equations in one space dimension:
\begin{align}\label{NLS}
 i\pa_t u_j+ \frac{1}{2m_j}\pa_x^2 u_j
 =
 F_j(u),
 \qquad
 (t,x)\in(0,\infty)\times\R
\end{align}
with the initial condition
\begin{align}\label{data}
 u_j(0,x) = u^\circ_j(x),\qquad x\in\R
\end{align}
for $j=1,\ldots,N$, where $i=\sqrt{-1}$,
$u=(u_j)_{1\le j\le N}$ is a $\C^N$-valued
unknown function of $(t,x)\in[0,\infty)\times\R$
and the masses $m_1,\ldots,m_N$ are positive constants.
Simply we assume that the nonlinear term $F=(F_j)_{1\le j\le N}:\C^N\to\C^N$ 
is a cubic homogeneous polynomial
in $\left(u, \conj u\right)$ 
with some complex coefficients, i.e.,
\[
 F_j(u)
 =
 \sum_{k,l,m=1}^N\,
 \sum_{\sigma\in\{+,-\}^3}
 C_{jklm}^\sigma
 u_k^{(\sigma_1)}
 u_l^{(\sigma_2)}
 u_m^{(\sigma_3)}
\]
with complex constants $C_{jklm}^{\sigma}$,
where $\sigma=(\sigma_1, \sigma_2, \sigma_3)$ and
$u_j^{(+)}=u_j$, $u_j^{(-)} = \conj{u_j}$ respectively.
Also we assume that the system satisfies
{\it gauge invariance}, i.e.,
\begin{align}\label{gauge}
 F_j\left(e^{im_1\theta}z_1, \ldots, e^{im_N\theta}z_N\right)
 =
 e^{im_j\theta}F_j\left(z_1,\ldots,z_N\right)
\end{align}
for each $j=1,\ldots,N$ 
and any $\theta\in\R$, $z=(z_j)_{1\le j\le N}\in\C^N$.
In the present paper, we are interested in large-time behavior of the
small amplitude solution for \eqref{NLS}--\eqref{data}.

Let us recall some previous results briefly.
There is a large body of literature discussing global existence
and large-time behavior of solutions for the single nonlinear Schr\"odinger equations in $n$-space dimensions of the form
\begin{align}\label{dNLS}
 i\pa_t u + \frac{1}{2}\Delta u = G(u),
 \qquad
 (t,x)\in\R\times\R^n,
\end{align}
where $\Delta$ is the Laplace operator in $x\in\R^n$
and $G(u)$ is a nonlinear term. 
We refer the readers to \cite{Cazenave03}
concerning the recent development on studies of \eqref{dNLS}.
Let us denote by $p_S(n)$ the Strauss exponent, which is defined by
$p_S(n)=\frac{n+2+\sqrt{n^2+12n+4}}{2n}$.
Strauss showed in \cite{Strauss81} that if the nonlinear term
$G(u)$ satisfies $|G'(u)|\le C|u|^{p-1}$ with $p>p_S(n)$,
then there exists a unique global solution for \eqref{dNLS}
with a suitable small initial data.
Note that for one-dimensional case, we have 
$p_S(1)=\frac{3+\sqrt{17}}{2}\approx 3.56$.
Now we concentrate our attention to 
the power-type nonlinearity, i.e. the case that $G(u)=\lambda |u|^{p-1}u$ with
$\lambda\in\C\setminus\{0\}$,
because it is a typical one satisfying the gauge invariance condition \eqref{gauge}.
In this case, when $p>1+2/n$, it is well-known that any solution 
$u(t)$ of \eqref{dNLS} behaves like a free solution 
as $t\to\infty$, if the data belongs to suitable weighted Sobolev spaces
(see e.g. \cite{Tsutsumi84}).
On the other hand, Barab \cite{Barab84} showed that there is no
asymptotically free solution for \eqref{dNLS}
if $1\le p\le 1+2/n$.
So we can see that cubic nonlinearities are critical for \eqref{dNLS} when $n=1$. 

Now we consider \eqref{dNLS} with $G(u)=\lambda|u|^2 u$,
$\lambda\in\C$ in one space dimension.
In this case, we can find an asymptotic profile of the solution to \eqref{dNLS}
in \cite{Ozawa91} for sufficiently small final data, if $\lambda\in\R$. 
We note that the asymptotic profile given there is just a phase-shifted free profile,
so the amplitude of the solution still behaves like a free solution.
Similar results can be found in \cite{Hayashi98}.
More precisely, Hayashi and Naumkin proved in \cite{Hayashi98} that the solution
decays like $O(t^{-1/2})$ in $L^\infty$ and behaves like a phase-shifted
free solution when $\lambda\in\R$, if the initial data is sufficiently small and
in suitable weighted Sobolev spaces.
On the other hand, according to the result by Shimomura \cite{Shimomura06}, 
the solution decays like $O(t^{-1/2}(\log t)^{-1/2})$ in $L^\infty$,
if $\im \lambda <0$ and the initial data is small enough.
Remember that the $L^\infty$-decay rate of the free evolution is $t^{-n/2}$
for $n$-dimensional cases.
Therefore this gain of additional logarithmic time-decay can
be read as a kind of the long-range effect.
His result was extended by \cite{Hayashi08}, which considers
the nonlinear terms including derivative types also.

Next we turn our attention to the case of systems,
i.e. \eqref{NLS} with $N\ge 2$ where $\pa_x^2$ is replaced by $\Delta$.
In this case, the problem becomes more complicated because
global existence and large-time behavior of the
solution are affected by the ratio of masses 
as well as the structure of nonlinearities.
Li found some structural conditions on the 
quadratic nonlinearities and the masses in \cite{Li12}
under which the solution to \eqref{NLS} 
exists globally and decays like a free solution in two space
dimension, if the data is small and belongs to a suitable weighted Sobolev space
(note that quadratic nonlinearities are critical in two space dimension).
And by a minor modification of the method there, 
we can obtain similar results for one-dimensional cases also.
Recently, Katayama, Li and Sunagawa considered
the quadratic two-dimensional NLS system 
\begin{align}\label{KLS}
 \left\{
 \begin{array}{l}
  i\pa_t u_1 + \dfrac{1}{2m_1}\Delta u_1 
  =\lambda_1|u_1|u_1+\mu_1 \conj{u_2} u_3,\\[1.8ex]
  i\pa_t u_2 + \dfrac{1}{2m_2}\Delta u_2 
  =\lambda_2|u_2|u_2+\mu_2 \conj{u_1} u_3,
  \qquad (t,x)\in (0,\infty)\times\R^2
  \\[1.8ex]
  i\pa_t u_3 + \dfrac{1}{2m_3}\Delta u_3 
  =\lambda_3|u_3|u_3 +\mu_3 u_1 u_2
 \end{array}
 \right.
\end{align}
in \cite{Katayama14} with complex coefficients 
on the nonlinearities under the mass relation $m_1+m_2=m_3$.
They showed if $\im \lambda_j <0$ for $j=1,2,3$ and
$\kappa_1 \mu_1 + \kappa_2 \mu_2 = \kappa_3 \conj{\mu_3}$
with some $\kappa_1, \kappa_2, \kappa_3 >0$, then
the solution of \eqref{KLS} decays like $O(t^{-1}(\log t)^{-1})$
in $L^\infty$ for sufficiently small data which belongs to certain
weighted Sobolev spaces. We refer the readers to
\cite{Li14-2} and the references cited therein 
for the recent progress on two-dimensional NLS systems
with critical nonlinearities.

This paper can be regarded as a one-dimensional version of the paper 
\cite{Katayama14} or a piece of extension of the paper \cite{Shimomura06}.
The aim of the present work is to introduce a structural condition 
of the cubic nonlinearities and the masses under which 
\eqref{NLS}--\eqref{data} admits a unique global solution
and it decays like $O(t^{-1/2}(\log t)^{-1/2})$
in $L^\infty$ as $t\to\infty$, 
if the initial data is small enough and belongs to suitable weighted
Sobolev spaces.

\section{Main Results} \label{sec_main}
In order to state our main results, we introduce some notations here.
We denote the usual Lebesgue space by $L^p(\R)$ equipped
with the norm
$\left\|\phi\right\|_{L^p} = \left(\int_\R|\phi(x)|^p\,dx\right)^{1/p}$
if $p\in [1,\infty)$ and $\left\|\phi\right\|_{L^\infty} = \sup_{x\in\R}|\phi(x)|$
if $p=\infty$.
The weighted Sobolev space is defined by
\[
 H_p^{s,q} (\R) = \left\{\phi = (\phi_1, \ldots, \phi_N)
 \in L^p(\R) : 
 \left\|\phi\right\|_{H_p^{s,q}}
 =
 \sum_{j=1}^N \left\|\phi_j\right\|_{H_p^{s,q}}<\infty\right\}
\]
with the norm
$\left\|\phi_j\right\|_{H_p^{s,q}} 
 = 
 \left\|\jb{x}^q \jb{i\pa_x}^s\phi_j\right\|_{L^p}$ 
for $s,q \in \R$ and $p\in[1,\infty]$,
where $\jb{\cdot} = \sqrt{1+\cdot\,^2}$.
For simplicity, we write $H^{s,q} = H_2^{s,q}$,
$H_p^s = H_p^{s,0}$
and the usual Sobolev space as $H^s = H^{s,0}$.
We define the Fourier transform $\hat\phi(\xi)$ of a function $\phi(x)$ by
\[
 (\F \phi)(\xi) = \hat \phi(\xi) 
 = 
 \frac{1}{\sqrt{2\pi}} \int_\R e^{-ix\xi}\phi(x)\,dx.
\]
Then the inverse Fourier transform is given by
\[
 (\IF\phi)(x)
 = 
 \frac{1}{\sqrt{2\pi}} \int_\R e^{ix\xi}\phi(\xi)\,d\xi.
\]
We denote by $y\cdot z$ the standard scalar product in $\C^N$
for $y,z \in \C^N$ and write $|z| = \sqrt{z\cdot z}$ as usual.
We can now formulate the main results.

\begin{thm}\label{thm_SDGE}
 Let $u^\circ\in H^{1,0}(\R)\cap H^{0,1}(\R)$ and
 $\left\|u^\circ\right\|_{H^{1,0}}+\left\|u^\circ\right\|_{H^{0,1}}=\eps$.
 Assume the condition \eqref{gauge} holds and suppose that
 there exists an $N\times N$ positive Hermitian matrix $A$ such that
 \begin{align}\label{condi_SDGE}
  \im\left(F(z)\cdot Az\right)
  \le
  0
 \end{align}
 for all $z\in\C^N$.
 Then there exists $\eps_0>0$ such that for all $\eps\in(0,\eps_0]$,
 the initial value problem \eqref{NLS}--\eqref{data} admits a unique
 global solution
 \[
 u(t)\in C^0\left([0,\infty);H^{1,0}(\R)\cap H^{0,1}(\R)\right)
 \]
 satisfying the time-decay estimate
 \[
  \left\|u(t)\right\|_{L^\infty}
  \le
  C(1+t)^{-1/2}
 \]
 for all $t\ge 0$.
\end{thm}

\begin{thm}\label{thm_decay}
 Suppose that the assumptions of Theorem~\ref{thm_SDGE} are fulfilled.
 Moreover, suppose that there exist an $N\times N$
 positive Hermitian matrix $A$ and constants $C_*,C^*>0$ such that 
 \begin{align}\label{condi_decay}
  -C_* |z|^4
  \le
  \im\left(F(z)\cdot Az\right)
  \le
  -C^* |z|^4
 \end{align}
 for all $z\in\C^N$.
 Then the global solution of \eqref{NLS}--\eqref{data},
 which is guaranteed by Theorem~\ref{thm_SDGE}, satisfies the time-decay
 estimate
 \[
  \left\|u(t)\right\|_{L^\infty}
  \le
  \frac{C(1+t)^{-1/2}}{\sqrt{\log(2+t)}}
 \]
 for all $t\ge 0$.
\end{thm}

\noindent
Here we give some examples satisfying the assumption 
of Theorem~\ref{thm_SDGE} or Theorem~\ref{thm_decay}
with suitable mass relations.\\

\noindent
{\bf Example 2.1.}
We consider the following two-component system
\begin{align}\label{ex1}
 \left\{
 \begin{array}{l}
  i\pa_t u_1 + \dfrac{1}{2m_1}\pa_x^2 u_1 
  =F_1(u)=\lambda_1|u|^2 u_1+\mu_1 \conj{u_1}^2 u_2,\\[2ex]
  i\pa_t u_2 + \dfrac{1}{2m_2}\pa_x^2 u_2 
  =F_2(u)=\lambda_2|u|^2 u_2 +\mu_2 u_1^3
 \end{array}
 \right.
\end{align}
in $(t,x)\in(0,\infty)\times\R$ 
under the mass relation $m_2 = 3m_1$,
where $\lambda_1, \lambda_2, \mu_1, \mu_2\in\C$.
Then we can see that the system \eqref{ex1} satisfies the gauge invariance
condition \eqref{gauge}.
Also we assume that the constants satisfy the following conditions:
\begin{align*}
 \im\lambda_j <0\;\;\text{for}\;\;j=1,2,\qquad
 \kappa_1\mu_1=\kappa_2\conj{\mu_2}
 \;\;\text{with some}\;\;\kappa_1,\kappa_2>0.
\end{align*}
Then we have 
\begin{align*}
 -C_*|z|^4
 \le
 \im\left(F(z)\cdot Az\right)
 =
 \sum_{j=1}^2 \kappa_j\left(\im\lambda_j\right)|z|^2|z_j|^2
 +
 \im\left(\kappa_1\mu_1\conj{z_1}^3 z_2
 +\kappa_2\mu_2 z_1^3\conj{z_2}\right)
 \le
 -C^*|z|^4
\end{align*}
for all $z\in\C^2$ with $A=\text{diag}\left(\kappa_1,\kappa_2\right)$.
Therefore we can conclude from Theorem~\ref{thm_SDGE} and
Theorem~\ref{thm_decay} that there exists a unique global solution
$u(t)\in C^0\left([0,\infty);H^{1,0}(\R)\cap H^{0,1}(\R)\right)$
to the initial value problem \eqref{ex1}--\eqref{data}
and the solution decays like $O\left(t^{-1/2}(\log t)^{-1/2}\right)$ in $L^\infty$
as $t\to\infty$, if $u^\circ\in H^{1,0}(\R)\cap H^{0,1}(\R)$ and it is 
sufficiently small.\\

\noindent
\textbf{Remark 2.1.}
 Since the Klein-Gordon equation is a relativisitic version of
 the Schr\"odinger equation, it is interesting to compare our results
 with a system of nonlinear Klein-Gordon equations.
 Here we consider the following two-component cubic nonlinear 
 Klein-Gordon system including dissipative nonlinearities
 \begin{align}\label{NLKG}
 \left\{
 \begin{array}{l}
  (\pa_t^2-\pa_x^2+m_1^2)u_1
  =\lambda_1|\pa_t u|^2 \pa_t u_1
    -(\pa_t u_1)^2 \pa_t u_2\\
  (\pa_t^2-\pa_x^2+m_2^2)u_2=
  \lambda_2|\pa_t u|^2 \pa_t u_2
    +(\pa_t u_1)^3
 \end{array}
 \right.
\end{align}
in $(t,x)\in(0,\infty)\times\R$ with the same mass relation 
$m_2 = 3m_1$ (often called {\it mass resonance} relation) as above, 
where $u=(u_1, u_2)$ is real-valued and
$\lambda_1, \lambda_2 <0$.
Then as pointed out in \cite{Kim15}, we can modify the proof of
\cite{Kim14} to see that \eqref{NLKG} admits a unique global solution $u(t)$
and it decays like $O\left(t^{-1/2}(\log t)^{-1/2}\right)$ in $L^\infty$
as $t\to\infty$, if the Cauchy data are sufficiently smooth, small
and compactly-supported.
We remark that the condition $\lambda_1, \lambda_2 <0$
reflects a dissipative character in this case, 
as the condition $\im \lambda_1, \im\lambda_2 <0$ implies a dissipative
property in \eqref{ex1}.\\

\noindent
We end this section by giving an example which satisfies the condition
\eqref{condi_SDGE} but violates \eqref{condi_decay}.\\

\noindent
{\bf Example 2.2.}
We consider the following four-component system 
\begin{align}\label{ex2}
 \left\{
 \begin{array}{l}
  i\pa_t u_1 + \dfrac{1}{2m_1}\pa_x^2 u_1 
  =F_1(u)=\mu_1 \conj{u_2} \conj{u_3} u_4,\\[1.8ex]
    i\pa_t u_2 + \dfrac{1}{2m_2}\pa_x^2 u_2
  =F_2(u)=\mu_2\conj{u_3} u_4 \conj{u_1},\\[1.8ex]
    i\pa_t u_3 + \dfrac{1}{2m_3}\pa_x^2 u_3 
  =F_3(u)=\mu_3u_4 \conj{u_1} \conj{u_2},\\[1.8ex]
    i\pa_t u_4 + \dfrac{1}{2m_4}\pa_x^2 u_4 
  =F_4(u)=\mu_4 u_1 u_2 u_3
 \end{array}
 \right.
\end{align}
in $(t,x)\in(0,\infty)\times\R$
under the mass relation $m_4 = m_1+m_2+m_3$, 
where $\mu_1, \ldots, \mu_4\in\C$.
Then we can see that the condition \eqref{gauge} holds and
$\im\left(F(z)\cdot Az\right)=0$ for all $z\in\C^4$ 
with $A=\text{diag}(\kappa_1,\kappa_2, \kappa_3, \kappa_4)$,
if there exist some positive constants $\kappa_1, \ldots, \kappa_4$
such that
$\kappa_1\mu_1 +\kappa_2\mu_2 +\kappa_3\mu_3 
= \kappa_4\conj{\mu_4}$.
Therefore in this case,
we can conclude from Theorem~\ref{thm_SDGE} that
the global solution of \eqref{ex2} 
decays like a free solution if the data is small enough.\\

The rest of this paper is organized as follows. 
In Section~\ref{sec_pre}, we compile some basic facts
concerning the free Schr\"odinger evolution group.
Section~\ref{sec_apriori} is devoted to obtain a suitable a priori estimate 
from which Theorem~\ref{thm_SDGE} follows immediately.
After that, we prove Theorem~\ref{thm_decay} in Section~\ref{sec_proof}
and discuss the optimality of the decay-rate of the solution.
In what follows, all non-negative constants will be denoted by $C$
which may vary from line to line unless otherwise specified.

\section{Preliminaries} \label{sec_pre}
In this section, we introduce some notations and useful estimates
which will be used in Section~\ref{sec_apriori} and Section~\ref{sec_proof}
for the proof of the main results.
In what follows, we denote
$\mathcal A\mathcal B = \left(\mathcal A_j \mathcal B_j\right)_{1\le j\le N}$
for $N$-dimensional column vectors $\mathcal A=(\mathcal A_j)_{1\le j\le N}$ 
and $\mathcal B=(\mathcal \B_j)_{1\le j\le N}$.
First we introduce the free Schr\"odinger evolution group 
$\U(t)=\left(\U_j(t)\right)_{1\le j\le N}$
defined by
\[
 \U_j(t)
 =
 e^{\frac{it}{2m_j}\pa_x^2}
 =
 \IF e^{-\frac{it}{2m_j}\xi^2}\F.
\]
It is well-known that $\U(t)$ is decomposed into
$\U(t)=\M(t)\D(t)\G\M(t)$, where the multiplication factor
$\M(t)=\left(\M_j(t)\right)_{1\le j\le N}$ is defined by
$\M_j(t)\phi(x)=\exp(\frac{im_j}{2t}x^2)\phi(x)$,
the Fourier-like transform $\G=(\G_j)_{1\le j\le N}$ 
(see e.g. \cite{Katayama14})
and the dilation operator $\D(t)$
are given by
\[
 (\G_j\phi)(\xi)
 =
 \sqrt{\frac{m_j}{i}}
 (\F\phi)(m_j \xi),
 \qquad
 \D(t)\phi(x)
 =
 \frac{1}{\sqrt{t}}\phi\left(x t^{-1}\right).
\]
Also we define $\W(t)=\G\M(t)\IG$ so that
$\U(t)=\M(t)\D(t)\W(t)\G$.
Then we have
\begin{align}\label{W-1Linfty}
 \left\|\left(\W(t)-1\right)\phi\right\|_{L^\infty}
 \le
 Ct^{-1/4}\left\|\phi\right\|_{H^1},
 \qquad
 \left\|\left(\W^{-1}(t)-1\right)\phi\right\|_{L^\infty}
 \le
 Ct^{-1/4}\left\|\phi\right\|_{H^1}
\end{align}
for $\phi=(\phi_1,\ldots,\phi_N)\in H^1(\R)$.
Indeed it is easy to check that the estimates \eqref{W-1Linfty} hold.
Since 
\[
 \left|\M_j(t)-1\right|
 =
 \left|e^{\frac{im_j}{2t}x^2}-1\right|
 =
 2\left|\sin\frac{m_j x^2}{4t}\right|
 \le
 C\frac{|x|^{2\beta}}{t^\beta}
\]
where $\beta\in[0,1]$, 
we find
\begin{align*}
 \left\|\left(\W(t)-1\right)\phi\right\|_{L^\infty}
 &=
 \left\|\G(\M(t)-1)\IG\phi\right\|_{L^\infty}
 \le
 C\left\|(\M(t)-1)\IG\phi\right\|_{L^1}\\
 &\le
 Ct^{-\beta}\left\|\jb{x}^{-\eta}\jb{x}^\eta|x|^{2\beta}\IF\phi\right\|_{L^1}
 \\
 &\le
 Ct^{-\beta}\left\|\jb{x}^{-\eta}\right\|_{L^2}
 \left\|\jb{x}^{2\beta+\eta}\IF\phi\right\|_{L^2}\\
 &\le
 Ct^{-\beta}\left\|\phi\right\|_{H^{2\beta+\eta}}
\end{align*}
holds for any $\eta>1/2$.
So by choosing $\beta=1/4$, we get the first estimate of \eqref{W-1Linfty}.
In view of the relation $\W^{-1}(t)-1 = -(\W(t)-1)\W^{-1}(t)$ and
$\left\|\W^{-1}(t)\phi\right\|_{H^1}\le C\left\|\phi\right\|_{H^1}$,
the second estimate of \eqref{W-1Linfty} follows immediately.

\section{A Priori Estimates} \label{sec_apriori}
The argument of this section is similar to those of the previous works,
for example
\cite{Katayama14}, \cite{Li12} and \cite{Shimomura06}.
Let $u(t)$ be the solution to \eqref{NLS}--\eqref{data} in $[0,T]$
and we define
\[
 \left\|u\right\|_{X_T}
 =
 \sup_{t\in[0,T]}
 \left(
 \jb{t}^{-\gamma}\big(
 \left\|u(t)\right\|_{H^1}+\left\|\U(-t)u(t)\right\|_{H^{0,1}}\big)
 +\jb{t}^{1/2}\left\|u(t)\right\|_{L^\infty}\right)
\]
where $0<\gamma\ll 1$ small.
\begin{lem}\label{lem_apriori}
 Under the assumption of Theorem~\ref{thm_SDGE}, 
 there exist $\eps_1>0$ and $C_0>0$ 
 such that $\left\|u\right\|_{X_T}\le\sqrt\eps$ implies
 $\left\|u\right\|_{X_T}\le C_0\eps$
for any $\eps\in(0,\eps_1]$. Here the constant $C_0$ does not
depend on $T$.
\end{lem}
Once this lemma is proved, we can obtain the global existence part of 
Theorem~\ref{thm_SDGE} in the following way:
By taking $\eps_0\in(0,\eps_1]$ so that $2C_0\sqrt{\eps_0}\le 1$,
we deduce that $\left\|u\right\|_{X_T}\le\sqrt\eps$
implies $\left\|u\right\|_{X_T}\le \sqrt\eps/2$
for any $\eps\in(0,\eps_0]$.
Then by the continuity argument, we have
$\left\|u\right\|_{X_T}\le C_0\eps$ as long as the solution exists.
So the local solution to \eqref{NLS}--\eqref{data}
can be extended to the global one.\\

From now on, we will prove Lemma~\ref{lem_apriori}. 
Since $i\pa_t \U_j(-t)=\U_j(-t)(i\pa_t+\frac{1}{2m_j}\pa_x^2)$,
taking $\U(-t)$ to the both sides of \eqref{NLS}, we get the following
integral equation
\[
 u(t)
 =
 \U(t)u^\circ
 -i\int_0^t\U(t-\tau)F(u(\tau))\,d\tau.
\]
Taking the $H^1$ norm, we obtain
\begin{align}\label{uH1}
 \left\|u(t)\right\|_{H^1}
 &\le
 \left\|u^\circ\right\|_{H^1}
 +\int_0^t\left\|F(u(\tau))\right\|_{H^1}\,d\tau\nonumber\\
 &\le
 C\eps+C\int_0^t\left\|u(\tau)\right\|_{L^\infty}^2
 \left\|u(\tau)\right\|_{H^1}\,d\tau\nonumber\\
 &\le
 C\eps+C\int_0^t\eps^{3/2}\jb{\tau}^{-1+\gamma}\,d\tau
 \le
 C\eps\jb{t}^\gamma,
\end{align}
where we used the assumption $\left\|u\right\|_{X_T}\le\sqrt\eps$.
Now we are going to estimate $\left\|\U(-t)u(t)\right\|_{H^{0,1}}$.
Since $\M(-t)F(u)=F(\M(-t)u)$ by \eqref{gauge}, 
using the relation 
\[
 \U_j(t)x\U_j(-t)=\M_j(t)itm_j^{-1}\pa_x\M_j(-t),
\]
we have
\begin{align}\label{xU-tsub}
 \left\|x\U(-t)F(u)\right\|_{L^2}
 &=
 \left\|\U(t)x\U(-t)F(u)\right\|_{L^2}
 \le
 Ct\left\|\pa_x\M(-t) F(u)\right\|_{L^2}\nonumber\\
 &=
 Ct\left\|\pa_x F(\M(-t)u)\right\|_{L^2}
 \le
 Ct\left\|\M(-t)u\right\|_{L^\infty}^2
 \left\|\pa_x \M(-t)u\right\|_{L^2}\nonumber\\
 &\le
 C\left\|u\right\|_{L^\infty}^2 \left\|x\U(-t)u\right\|_{L^2}.
\end{align}
Therefore by the similar way as above, we obtain
\begin{align}\label{U(-t)uH01}
 \left\|\U(-t)u(t)\right\|_{H^{0,1}}
 &\le
 \left\|u^\circ\right\|_{H^{0,1}}
 +\int_0^t
 \left\|\U(-\tau)F(u(\tau))\right\|_{H^{0,1}}\,d\tau\nonumber\\
 &\le
 C\eps+C\int_0^t\left\|u(\tau)\right\|_{L^\infty}^2
 \left\|\U(-\tau)u(\tau)\right\|_{H^{0,1}}\,d\tau\nonumber\\
 &\le
 C\eps\jb{t}^\gamma,
\end{align}
where we used \eqref{xU-tsub} for the second inequality. 
Next, we consider the term $\left\|u(t)\right\|_{L^\infty}$.
If $t\le 1$, the standard Sobolev embedding 
and \eqref{uH1} suffice to obtain
\begin{align}\label{uLinftysmallt}
 \jb{t}^{1/2+\gamma}\left\|u(t)\right\|_{L^\infty}
 \le
 2^{1/4+\gamma/2}\left\|u(t)\right\|_{H^1}
 \le
 C\eps\jb{t}^\gamma.
\end{align}
So from now on we consider the case that $t\ge 1$.
We define the new function $\alpha=(\alpha_j)_{1\le j\le N}$ by
\[
 \alpha(t,\xi)
 =
 \G\left(\U(-t)u(t)\right)(\xi).
\]
Then from the decomposition of the free Schr\"odinger evolution group
and \eqref{gauge}, we have
\begin{align}\label{ipatalpha}
 i\pa_t\alpha(t,\xi)
 &=
 \G\U(-t)F(u(t))(\xi)\nonumber\\
 &=
 \W^{-1}(t)\D^{-1}(t)\M(-t)
 F(\M(t)\D(t)\W(t)\alpha(t,\xi))\nonumber\\
 &=
 t^{-1}\W^{-1}(t)F(\W(t)\alpha(t,\xi))\nonumber\\
 &=
 t^{-1}F(\alpha(t,\xi))+R(t,\xi),
\end{align}
where
\[
 R(t,\xi)
 =
 t^{-1}\W^{-1}(t)F(\W(t)\alpha(t,\xi))-t^{-1}F(\alpha(t,\xi)).
\]
As we shall see below, 
$R$ can be regarded as a remainder because it decays 
strictly faster than $O(t^{-1})$ in $L^\infty$, while the first term of 
the right-hand side of \eqref{ipatalpha} plays a role as a main term.
Since we have by the estimate \eqref{W-1Linfty}, the Sobolev embedding
and \eqref{U(-t)uH01},
\begin{align*}
 t^{-1}\left\|\left(\W^{-1}(t)-1\right)F(\W(t)\alpha(t))\right\|_{L^\infty}
 &\le
 Ct^{-5/4}\left\|F(\W(t)\alpha(t))\right\|_{H^1}
 \le
 Ct^{-5/4}\left\|\W(t)\alpha(t)\right\|_{H^1}^3\\
 &\le
 Ct^{-5/4}\left\|\alpha(t)\right\|_{H^1}^3
 \le
 Ct^{-5/4}\left\|\G\U(-t)u(t)\right\|_{H^1}^3\\
 &\le
 Ct^{-5/4}\left\|\U(-t)u(t)\right\|_{H^{0,1}}^3\\
 &\le
 C\eps^3 t^{-5/4+3\gamma}
\end{align*}
and similarly
\begin{align*}
 &t^{-1}\left\|F(\W(t)\alpha(t))-F(\alpha(t))\right\|_{L^\infty}\\
 &\qquad\le
 Ct^{-1}\left\|(\W(t)-1)\alpha(t)\right\|_{L^\infty}
 \left(\left\|(\W(t)-1)\alpha(t)\right\|_{L^\infty}^2
 +\left\|\W(t)\alpha(t)\right\|_{L^\infty}\left\|\alpha(t)\right\|_{L^\infty}\right)\\
 &\qquad\le
 Ct^{-5/4}\left\|\alpha(t)\right\|_{H^1}
 \left(t^{-1/2}\left\|\alpha(t)\right\|_{H^1}^2
 +\left\|\alpha(t)\right\|_{H^1}^2\right)\\
 &\qquad\le
 C\eps^3 t^{-5/4+3\gamma},
\end{align*}
we deduce that the remainder satisfies the estimate
\begin{align}\label{remainder}
 \left\|R(t)\right\|_{L^\infty}
 &\le
 t^{-1}
 \left(\left\|\left(\W^{-1}(t)-1\right)F(\W(t)\alpha(t))\right\|_{L^\infty}
 +
 \left\|F(\W(t)\alpha(t))-F(\alpha(t))\right\|_{L^\infty}
 \right)\nonumber\\
 &\le
 C\eps^3 t^{-5/4+3\gamma}
\end{align}
for $t\ge 1$.
Here we note that
\begin{align}\label{matA}
 |y\cdot Az| \le (y\cdot Ay)^{1/2}(z\cdot Az)^{1/2},
 \qquad
 \lambda_* |z|^2 \le z\cdot Az \le \lambda^* |z|^2
\end{align}
hold for $y,z\in\C^N$, where the matrix $A$ is in Theorem~\ref{thm_SDGE}
and $\lambda^*$ (resp. $\lambda_*$) is the largest (resp. smallest) 
eigenvalue of $A$.
Therefore it follows from \eqref{ipatalpha}, \eqref{condi_SDGE},
\eqref{matA} and \eqref{remainder} that
\begin{align}\label{patalphaAalpha}
 \pa_t\left(\alpha(t,\xi)\cdot A\alpha(t,\xi)\right)
 &=
 2\im\left(i\pa_t\alpha(t,\xi)\cdot A\alpha(t,\xi)\right)\nonumber\\
 &=
 2t^{-1}\im\left(F(\alpha(t,\xi))\cdot A\alpha(t,\xi)\right)
 +2\im\left(R(t,\xi)\cdot A\alpha(t,\xi)\right)\nonumber\\
 &\le
 Ct^{-5/4}\left|t^{5/4} R(t,\xi)\cdot A\alpha(t,\xi)\right|\nonumber\\
 &\le
 Ct^{-5/4}\left(\alpha(t,\xi)\cdot A\alpha(t,\xi)
 +t^{5/2}R(t,\xi)\cdot AR(t,\xi)\right)\nonumber\\
 &\le
 Ct^{-5/4}\alpha(t,\xi)\cdot A\alpha(t,\xi)
 +C\eps^6 t^{-5/4+6\gamma}.
\end{align}
Noting that
\[
 \left\|\alpha(1)\right\|_{L^\infty}
 =
 \left\|\G\U(-1)u(1)\right\|_{L^\infty}
 \le
 C\left\|\G\U(-1)u(1)\right\|_{H^1}
 \le
 C\left\|\U(-1)u(1)\right\|_{H^{0,1}}
 \le
 C\eps,
\]
integrating \eqref{patalphaAalpha} with respect to time lead to
\[
 \alpha(t,\xi)\cdot A\alpha(t,\xi)
 \le
 C\eps^2 + \int_1^t \tau^{-5/4} \alpha(t,\xi)\cdot A\alpha(t,\xi)\,d\tau.
\]
Thus the Gronwall lemma and \eqref{matA} yield
\begin{align}\label{alphaLinfty}
 \left\|\alpha(t)\right\|_{L^\infty}
 \le
 C\eps
\end{align}
for $t\ge 1$.
Hence with the estimates \eqref{alphaLinfty}, \eqref{U(-t)uH01}
and the inequality
\begin{align}\label{FM-1UuLinfty}
 \left\|\G(\M(t)-1)\U(-t)u(t)\right\|_{L^\infty}
 &=
 \left\|(\W(t)-1)\G\U(-t)u(t)\right\|_{L^\infty}\nonumber\\
 &\le
 Ct^{-1/4}\left\|\G\U(-t)u(t)\right\|_{H^1},
\end{align}
we finally obtain
\begin{align}\label{uLinfty}
 \left\|u(t)\right\|_{L^\infty}
 &=
 \left\|\M(t)\D(t)\G\M(t)\U(-t) u(t)\right\|_{L^\infty}\nonumber\\
 &\le
 Ct^{-1/2}\left\|\G\M(t)\U(-t)u(t)\right\|_{L^\infty}\nonumber\\
 &\le
 Ct^{-1/2}\big(
 \left\|\G\U(-t)u(t)\right\|_{L^\infty}
 +\left\|\G(\M(t)-1)\U(-t)u(t)\right\|_{L^\infty}\big)\nonumber\\
 &\le
 Ct^{-1/2}\left(
 \left\|\alpha(t)\right\|_{L^\infty}
 +t^{-1/4}\left\|\U(-t)u(t)\right\|_{H^{0,1}}\right)\nonumber\\
 &\le
 C\eps t^{-1/2}
\end{align}
for $t\ge 1$.
By \eqref{uH1}, \eqref{U(-t)uH01}, \eqref{uLinftysmallt} and 
\eqref{uLinfty}, we arrive at Lemma~\ref{lem_apriori}
and the $L^\infty$-decay estimate in Theorem~\ref{thm_SDGE}
follows immediately.\qed

\section{Proof of Theorem~\ref{thm_decay}} \label{sec_proof}
Now we are in a position to prove Theorem~\ref{thm_decay}.
Note that the similar arguments of this section are also used in
the previous works \cite{Katayama14}, \cite{Katayama15} and \cite{Kim14}.
If $t\le 2$ then by \eqref{uH1} we have
\[
 \jb{t}^{1/2+\gamma}\sqrt{\log (2+t)}
 \left\|u(t)\right\|_{L^\infty}
 \le
 C\left\|u(t)\right\|_{H^1}
 \le
 C\eps\jb{t}^\gamma,
\] 
so we only consider the case that $t\ge 2$.
First we note that
\[
 \pa_t\left((\log t)^2\left(\alpha(t,\xi)\cdot A\alpha(t,\xi)\right)\right)
 =
 (\log t)^2 \pa_t(\alpha(t,\xi)\cdot A\alpha(t,\xi))
 +
 \frac{2\log t}{t}\alpha(t,\xi)\cdot A\alpha(t,\xi).
\]
Similarly to \eqref{patalphaAalpha}, we have
\begin{align*}
 \pa_t(\alpha(t,\xi)\cdot A\alpha(t,\xi))
 &=
 2t^{-1}\im\left(F(\alpha(t,\xi))\cdot A\alpha(t,\xi)\right)
 +2\im\left(R(t,\xi)\cdot A\alpha(t,\xi)\right)\\
 &\le
 -2C^* t^{-1}|\alpha(t,\xi)|^4
 +C|R(t,\xi)||\alpha(t,\xi)|\\
 &\le
 -2C^* t^{-1}|\alpha(t,\xi)|^4
 +C\eps^4 t^{-5/4+3\gamma},
\end{align*}
where we used
\eqref{remainder}, \eqref{matA}, \eqref{alphaLinfty} and the assumption \eqref{condi_decay}
with the constant $C^*$ appearing in Theorem~\ref{thm_decay}.
Also we have
\begin{align*}
 \alpha(t,\xi)\cdot A\alpha(t,\xi)
 &\le
 \lambda^*|\alpha(t,\xi)|^2
 =
 \frac{\lambda^*}{\sqrt{2C^*\log t}}
 \sqrt{2C^*\log t}|\alpha(t,\xi)|^2\\
 &\le
 \frac{(\lambda^*)^2}{4C^*\log t}
 +C^* \log t |\alpha(t,\xi)|^4
\end{align*}
by the Young inequality.
Piecing them together, we obtain
\[
 \pa_t\left((\log t)^2\left(\alpha(t,\xi)\cdot A\alpha(t,\xi)\right)\right)
 \le
 Ct^{-1}
 +C\eps^4 (\log t)^2 t^{-5/4+3\gamma}.
\]
Integrating with respect to time, we get
\begin{align*}
 (\log t)^2 \left(\alpha(t,\xi)\cdot A\alpha(t,\xi)\right)
 \le
 C\eps^2 + C\int_2^t \left(\tau^{-1}
 +\frac{(\log \tau)^2}{\tau^{5/4-3\gamma}}\right)\,d\tau
 \le
 C\log t
\end{align*}
for $t\ge 2$. Thus \eqref{matA} yields
\[
 \left\|\alpha(t)\right\|_{L^\infty}
 \le
 C(\log t)^{-1/2}.
\]
Therefore by the same arguments as in \eqref{uLinfty},
we arrive at
\begin{align*}
 \left\| u(t)\right\|_{L^\infty}
 \le
 Ct^{-1/2}\left(\left\|\alpha(t)\right\|_{L^\infty}
 +t^{-1/4}\left\|\U(-t)u(t)\right\|_{H^{0,1}}\right)
 \le
 Ct^{-1/2}(\log t)^{-1/2},
\end{align*}
which proves Theorem~\ref{thm_decay}.\qed\\

Finally, we discuss the optimality of the decay rate $O(t^{-1/2}(\log t)^{-1/2})$.
We put $u^\circ (x)= \delta v^\circ (x)
(\not\equiv 0)\in H^{1,0}(\R)\cap H^{0,1}(\R)$
with $\delta >0$
(note that $\eps = \left\|u^\circ\right\|_{H^{1,0}}
+\left\|u^\circ\right\|_{H^{0,1}}\le C\delta$).
Here we will show that the solution does not decay strictly faster than 
$t^{-1/2}(\log t)^{-1/2}$ as $t\to \infty$, if $\delta$ is sufficiently small.
Suppose that
\begin{align}\label{optimal}
 \lim_{t\to\infty} (t\log t)^{1/2} \left\|u(t)\right\|_{L^\infty}
 =
 0
\end{align}
holds.
By \eqref{FM-1UuLinfty} and \eqref{U(-t)uH01}, we have
\begin{align*}
 \left\|\alpha(t)\right\|_{L^\infty}
 &=
 \left\|\G\U(-t)u(t)\right\|_{L^\infty}\\
 &\le
 \left\|\G\M(t)\U(-t)u(t)\right\|_{L^\infty}
  +\left\|\G(\M(t)-1)\U(-t)u(t)\right\|_{L^\infty}\\
 &\le
 Ct^{1/2}\left\|\M(t)\D(t)\G\M(t)\U(-t)u(t)\right\|_{L^\infty}
 +Ct^{-1/4}\left\|\U(-t)u(t)\right\|_{H^{0,1}}\\
 &\le
 Ct^{1/2}\left\|u(t)\right\|_{L^\infty}
 +C\delta t^{-1/4+\gamma}.
\end{align*}
Thus from \eqref{optimal}, we get
\begin{align}\label{limit}
 (\log t)^{1/2}|\alpha(t,\xi)|
 \le
 C(t\log t)^{1/2}\left\|u(t)\right\|_{L^\infty}
 +C\delta t^{-1/4+\gamma}(\log t)^{1/2}
 \to
 0
\end{align}
as $t\to\infty$ uniformly with respect to $\xi\in\R$.
Hence if $\delta$ is sufficiently small, we have
\begin{align}\label{2Clambda}
 \sqrt{\frac{2C_*}{\lambda_*}}(\log t)^{1/2}|\alpha(t,\xi)|\le 1
\end{align}
for all $t\ge 2$ and $\xi\in\R$, where $C_*$ and $\lambda_*$
are the constants appearing in \eqref{condi_decay} and \eqref{matA} respectively.
Therefore as in \eqref{patalphaAalpha},
it follows from \eqref{condi_decay}, \eqref{matA}, 
\eqref{remainder}, \eqref{alphaLinfty} and \eqref{2Clambda} that
\begin{align*}
 &\pa_t \left( (\log t)\alpha(t)\cdot A\alpha(t)\right)
 =
 (\log t)\pa_t (\alpha(t)\cdot A\alpha(t))
 + t^{-1}\alpha(t)\cdot A\alpha(t)\\
 &\qquad=
 2(\log t)\left(t^{-1}\im\left(F(\alpha(t))\cdot A\alpha(t)\right)
 +\im \left(R(t)\cdot A\alpha(t)\right)\right)
 +t^{-1}\alpha(t)\cdot A\alpha(t)\\
 &\qquad\ge
 -2C_*t^{-1}(\log t)|\alpha(t)|^4
 +\lambda_*t^{-1}|\alpha(t)|^2
 -2(\log t)|R(t)\cdot A\alpha(t)|\\
 &\qquad\ge
 \lambda_* t^{-1}|\alpha(t)|^2
 \left(1-\frac{2C_*\log t }{\lambda_*}|\alpha(t)|^2\right)
 -C\delta^4 t^{-5/4+3\gamma}\log t\\
 &\qquad\ge
 -C\delta^4 t^{-5/4+3\gamma}\log t,
\end{align*}
which yields
\begin{align*}
 (\log t)\alpha(t)\cdot A\alpha(t)
 \ge
 (\log 2) \alpha(2)\cdot A\alpha(2)
 -C\delta^4\int_2^t\frac{\log \tau}{\tau^{5/4-3\gamma}}\,d\tau
 \ge
 C\delta^2 - C'\delta^4>0
\end{align*}
for sufficiently small $\delta$ with some positive constants $C$ and $C'$.
This contradicts \eqref{limit}.


\end{document}